# Cool WENO schemes

I. Cravero[a], G. Puppo[b,*], M. Semplice[a], G. Visconti[c]

[a]*Università di Torino. Via C. Alberto, 10. Torino, Italy.*
[b]*Università dell'Insubria. Via Valleggio, 10. Como, Italy.*
[c]*RWTH Aachen University. Templergraben 55, 52062 Aachen, Germany*

**Abstract**

This work is dedicated to the development and comparison of WENO-type reconstructions for hyperbolic systems of balance laws. We are particularly interested in high order shock capturing non-oscillatory schemes with uniform accuracy within each cell and low spurious effects. We need therefore to develop a tool to measure the artifacts introduced by a numerical scheme. To this end, we study the deformation of a single Fourier mode and introduce the notion of distorsive errors, which measure the amplitude of the spurious modes created by a discrete derivative operator. Further we refine this notion with the idea of temperature, in which the amplitude of the spurious modes is weighted with its distance in frequency space from the exact mode. Following this approach linear schemes have zero temperature, but to prevent oscillations it is necessary to introduce nonlinearities in the scheme, thus increasing their temperature. However it is important to heat the linear scheme just enough to prevent spurious oscillations. With several tests we show that the newly introduced CWENOZ schemes are cooler than other existing WENO-type operators, while maintaining good non-oscillatory properties.

## 1. Introduction

We consider hyperbolic systems of balance laws of the form

$$\partial_t u + \nabla \cdot f(u) = S(u), \qquad (1)$$

where $f$ is the flux function, and $S(u)$ is the source term. A very popular technique to integrate numerically systems of conservation and balance laws is the finite volume method. Here one splits the computational domain into the union of finite cells $\Omega_j$ and evolves the cell averages $\overline{u}_j(t)$ of the solution. In particular, projecting the equation on the grid and integrating by parts, one has to integrate the system of ODEs

$$\frac{d\overline{u}_j}{dt} = -\frac{1}{|\Omega_j|} \int_{\partial\Omega_j} f(u(t,s)) \cdot n(s) ds + \frac{1}{|\Omega_j|} \int_{\Omega_j} S(u(x)) dx \qquad (2)$$

where $\partial\Omega_j$ is the boundary of the cell $\Omega_j$, and $n$ is the outward unit normal to $\partial\Omega_j$. When a numerical scheme is introduced, the above equation must be closed choosing a technique to compute the point values of the unknown function $u$ from the cell averages to selected point values. More in detail, the point values of $u$ are needed at the nodes of a quadrature formula along $\partial\Omega_j$ for evaluating the contour integral of $f$, and at the nodes of a quadrature formula within $\Omega_j$, to approximate the integral of the source $S$.

*Corresponding author
Email addresses: `isabella.cravero@unito.it` (I. Cravero),
`gabriella.puppo@uninsubria.it` (G. Puppo), `matteo.semplice@unito.it` (M. Semplice),
`visconti@igpm.rwth-aachen.de` (G. Visconti)

A class of high order numerical methods to reconstruct these point values is WENO (Weighted Essentially Non Oscillatory), introduced in a general framework in [JS96], with the successive reviews [Shu98, Shu09]. These schemes are designed to optimize accuracy, when the solution is smooth, while decreasing accuracy in the presence of singularities in the data to avoid spurious oscillations. The local computational grids of WENO schemes are composed of $r$ overlapping substencils of $r$ cells, forming a larger stencil with $2r - 1$ cells. The method exploits the whole stencil with $2r - 1$ cells, when the data are smooth, achieving its optimal accuracy $p = 2r - 1$. If a discontinuity occurs within the large stencil, the scheme automatically avoids the discontinuity selecting one of the $r$ substencils, and as a consequence decreasing accuracy down to order $r$. The nonlinear coefficients of WENO's convex combination are based on lower order local smoothness indicators that are a scaled sum of the square $L^2$ norms of all derivatives of local interpolating polynomials, [JS96]. At smooth parts of the solution, all the smoothness indicators are roughly equal and the WENO weighted combination reproduces the central scheme of the polynomial of maximum degree, i.e. $2r - 2$. An essentially zero weight is assigned to those lower degree polynomials whose substencils contain high gradients or shocks, aiming at an essentially non-oscillatory solution close to discontinuities.

It was noticed that, close to critical points, the accuracy of WENO schemes may be non-optimal since the algorithm produces in this case non linear weights that are too far apart from the optimal values. In an attempt to alleviate this problem, first a remapping technique called WENOM was introduced in [HAP05] and then a new algorithm for the computation of the nonlinear weights was published in [BCCD08] and later extended in [CCD11, DB13]. This latter technique, called WENOZ, relies on an extra global higher order smoothness indicator that enters in the definition of the nonlinear weights in order to drive their values closer to the optimal ones in the case of smooth data. Since this extra indicator is simply defined as a linear combination of the Jiang-Shu indicators already computed, the extra cost is negligible. A study of the performance of the WENOZ reconstruction then revealed that this idea not only solves the issue of the accuracy at critical points, but also yields a more accurate reconstruction than its WENO counterpart, [CCD11, DB13].

However, WENOZ shares with WENO the feature of replicating the optimal polynomial only at one point at a time. If the reconstruction is needed at several points, as in the quadratures required by the integration of (2), then several reconstruction steps must be computed, each time with different weights. CWENO reconstructions were introduced in [LPR99] in order to reconstruct the point value at the center of the cell with third order accuracy in a staggered central scheme. Later, this technique has been extended to non-staggered finite volume schemes, to fifth [Cap08] and to any order [CPSV16] and exploited for general conservation and balance laws. See also [BGS16]. At a difference with WENO, CWENO constructs a polynomial $P_{\text{rec}}$ of degree $2r-2$, which is defined globally in the reconstruction cell. $P_{\text{rec}}$ is computed as a convex linear combination of the same $r$ polynomials of degree $r-1$ that are considered by WENO and of an extra polynomial $P_0$ of degree $2r-2$. The same mechanism used by WENO for computing the nonlinear weights of the convex combination is employed also by CWENO in order to guarantee that $P_{\text{rec}}$ has optimal accuracy $2r-1$ on smooth flows and is essentially non-oscillatory close to discontinuities. Due to its flexibility and uniform accuracy, this reconstruction technique has found many applications. For example in balance laws it allows to compute point value reconstructions at points in the interior of the cell (see [CPSV16]) without incurring the problem of non-existence of the WENO optimal weights described in [QS02, p. 194], without resorting to the treatment of non-positive optimal weights of [SHS02] and without the complication of having optimal weights that depend of the local mesh geometry as in [WFS08, PS16]. The freedom to choose the values of the CWENO optimal weights has been appreciated for schemes on locally refined cartesian grids (compare [CRS16] with e.g. [HS99, DK07]) and for moving mesh schemes like the ALE-ADER schemes on triangular and tetrahedral meshes [BSDR16].

Given the increased resolution and accuracy of WENOZ schemes when compared to their WENO counterparts (see [DB13]), in this paper we introduce a variant of the CWENO reconstruction, that we call CWENOZ and that makes use of the definition of the non linear weights defined by [DB13]. This is the subject of Section 2.

In order to compare different schemes, one may consider tables containing data on convergence histories on some classical problems, but they often depend heavily on the particular test problem chosen. Another standard criterium involves the concepts of numerical diffusion and numerical dispersion. A linear scheme applied to the linear advection equation propagates each single Fourier mode, modifying only the amplitude of each wave number and its propagation speed, but without generating spurious modes. This is measured with the concepts of numerical diffusion and numerical dispersion, and it can be exactly computed with the classical von Neumann analysis. [Pir06] proposed to extend the concepts of diffusion and dispersion to nonlinear schemes, by studying the real and imaginary part of the diagonal of the time-advancement operator in frequency space. In Section 3, we thus study the diffusion and dispersion of the WENO, CWENO and CWENOZ reconstructions, following the ideas of [Pir06], but applying them to the spatial operator only, in order to avoid the influence of the time advancement scheme.

The above generalizes the classical von Neumann analysis to non linear schemes, but it is still based on ideas which are derived from linear concepts, while essentially non-oscillatory schemes are highly non linear. In fact, when the discrete WENO, CWENO or CWENOZ derivative is applied to a single Fourier mode, spurious wave numbers are created. These spurious modes are then propagated during the time evolution, interacting in a highly non-trivial way. We thus think that it is important to study the size of these spurious modes. In order to evaluate this effect, in the second part of Section 3, we develop the idea of a distortion error that measures the size of the off-diagonal elements of the discrete derivative in frequency space. Moreover, we also consider a new concept, which we call *temperature*, that takes into account also the distance of the spuriously created modes from the original one in Fourier space. Both definitions incorporate a scaling factor, such that the computed quantities do not depend on the number of Fourier modes used in the computation and thus become a *signature* of the discrete derivative operator. As far as we know, this is the first attempt to measure the non linear distortion introduced by high order essentially non-oscillatory schemes.

We remark that a linear scheme is diagonal in frequency space and thus has zero distorsion and zero temperature. However, it is known that linear high order schemes are oscillatory close to singularities, which may arise in a finite time in the evolution of conservation laws. Thus a non-oscillatory scheme must have non zero temperature, i.e. it cannot be *cold*. On the other hand a scheme with high temperature is a scheme that may produce strong spurious signals and consequently large errors. One thus aims at *cool* schemes, in the sense that their temperature should be non-zero in order to be stable, but not too high in order to prevent excessive distorsion. We end Section 3 with a comparison of the distorsion and temperature of WENO, CWENO and CWENOZ schemes.

The paper is complemented in Section 4 by numerical experiments that confirm the results of the previous section and some conclusions are drawn in Section 5.

## 2. WENO-like reconstructions and the new CWENOZ

We start this section by reviewing essentially non-oscillatory reconstructions of WENO type. We recall that the mission of a reconstruction algorithm is to match the contrasting requirements of high accuracy one the one side and to avoid spurious oscillations on the other one. The WENO schemes of [JS96] have been a particularly successfull answer to this needs and have prompted the construction of variants designed to overcome some of their shortcomings. After the definition of the WENO reconstruction, we recall the CWENO reconstruction which ensures uniform accuracy within the cell and WENOZ which is devised to optimize the choice of the nonlinear WENO weights. Finally we blend CWENO and WENOZ, defining the new CWENOZ and conclude the section by proving its accuracy properties.

For an integer $r > 1$, let us consider reconstructions with a stencil of $2r - 1$ cells, centered on a cell $\Omega_0$, i.e. the cells with indices in $\mathbb{S}_0 = \{-r+1, \ldots, r-1\}$. Let $P_{\text{opt}}$ be the polynomial of degree $2r - 2$ that interpolates exactly all the cell averages in the stencil $\mathbb{S}_0$. This is clearly the most accurate reconstruction polynomial, if the data in $\mathbb{S}_0$ come from a smooth function. Typically essentially non-oscillatory recontructions, break $\mathbb{S}_0$ into $r$ smaller stencils $\mathbb{S}_k = \{-r+k, \ldots, k-1\}$ for $k = 1, \ldots, r$ and to consider the interpolating polynomials

$$P_k \in \mathbb{P}^{r-1} \text{ s.t. } \forall j \in \mathbb{S}_k : \int_{\Omega_j} P_k(s) \mathrm{d}s = h\overline{u}_j.$$

The WENO reconstruction at a point $\hat{x} \in \Omega_0$ consists in finding an optimal set of linear coefficients $d_k \in (0, 1)$ such that

$$P_{\text{opt}}(\hat{x}) = \sum_{k=1}^{r} d_k(\hat{x}) P_k(\hat{x})$$

and computing the reconstructed value

$$\text{WENO}(2r - 1; \hat{x}) = \sum_{k=1}^{r} \omega_k(\hat{x}) P_k(\hat{x}),$$

where the $\omega_k$ are the so-called nonlinear weights that are computed with the help of oscillation indicators $I[P_k]$ in such a way that $\omega_k \simeq d_k$ for smooth data, but $\omega_k \simeq 0$ if a discontinuity is present in stencil $\mathbb{S}_k$. See [JS96, Shu98, Shu09].

One of the drawbacks of the WENO reconstruction is the difficulty of guaranteeing the existence and positivity of a set of linear weights for an arbitrary point in the reconstruction cell. For example [QS02, p. 194] shows that, on uniform meshes, for even $r$, optimal weights for the cell center do not exist and, for odd $r$, they exist but are not in $(0, 1)$. There is a technique to treat negative weights [SHS02], but it requires to compute two different reconstructions per point.

In [LPR99], a new third order reconstruction technique was introduced, so that point values at cell centre could be computed. In fact, this technique, called Central WENO (CWENO), yields reconstructed values that are uniformly accurate in every point of the cell. It employs directly also an additional polynomial $P_0$ of degree $2r - 2$ in the linear combination. The technique has been later extended to treat higher order cases[CPSV16], non-uniform meshes and two and three-dimensional meshes of quads [SCR16] and simplices[BSDR16].

We recall here the Definition of the CWENO reconstruction, as given in [CPSV16], specialized to the present case.

**Definition 1** (CWENO). *For $r > 1$, consider the stencils $\mathbb{S}_0$ and the substencils $\mathbb{S}_1, \cdots, \mathbb{S}_r$ defined above, the polynomial $P_{\text{opt}}$ interpolating all the data in $\mathbb{S}_0$ and the polynomials $P_k$ interpolating the data in $\mathbb{S}_k$. Let also $d_0, d_1, \cdots, d_r$ be a set of positive, real coefficients $d_k \in (0, 1)$, s.t. $\sum_{k=0}^{r} d_k = 1$. The CWENO operator computes a reconstruction polynomial*

$$P_{\text{rec}} = \text{CWENO}(P_{\text{opt}}, P_1, \ldots, P_r) \in \mathbb{P}^{2r-2}$$

*as follows:*

1. *first, introduce the polynomial $P_0$ defined as*

$$P_0(x) = \frac{1}{d_0}\left(P_{\mathsf{opt}}(x) - \sum_{k=1}^{r} d_k P_k(x)\right) \in \mathbb{P}^{2r-2}; \qquad (3)$$

2. *then compute the nonlinear coefficients $\omega_k$ from the linear coefficients $d_k$ as*

$$\alpha_k = \frac{d_k}{(I[P_k] + \epsilon)^t}, \qquad \omega_k = \frac{\alpha_k}{\sum_{i=0}^{r} \alpha_i}, \qquad (4)$$

*where $I[P_k]$ denotes a suitable regularity indicator (e.g. the Jiang-Shu ones of eq. (9)) evaluated on the polynomial $P_k$, $\epsilon$ is a small positive quantity and $t \geq 2$;*

3. *and finally define*

$$P_{\mathsf{rec}}(x) = \sum_{k=0}^{r} \omega_k P_k(x) \in \mathbb{P}^{2r-2}. \qquad (5)$$

**Remark 1.** Note that the polynomial $P_{\mathsf{rec}}$ computed by CWENO can then be evaluated at reconstruction points and will have uniform accuracy for every $\hat{x} \in \Omega_0$. In particular, it will be $(2r-1)$-th order accurate on smooth data. This is particularly advantageous when source terms of balance laws are integrated.

Another advantage of CWENO over the classical WENO reconstruction is the invariance of the optimal linear coefficients $d_k$ on the reconstruction point. This feature, on top of ensuring the above-mentioned uniform accuracy, also allows to compute the nonlinear weights with equation (4) only once per cell and not once per reconstruction point.

Moreover, the optimal linear coefficients $d_k$ of CWENO are not bound by accuracy requirements as those of WENO. This means that they can be chosen also independently of the local mesh geometry (relative position and size) in non-uniform grids in one or more space dimension (compare e.g. [PS16, DK07]).

However CWENO, as well as WENO, requires one to choose $\varepsilon \asymp h^2$ or $\varepsilon \asymp h$ in order to be convergent with optimal order also on local extrema. In the WENO framework, Borges, Carmona, Costa and Don introduced the WENOZ reconstruction technique at order 5 in [BCCD08], which later was extended to arbitrary order in [CCD11]. This technique relies on the definition of an extra smoothness indicator $\tau$. The computation of $\tau$ is cheaper than the mapping technique of [HAP05], since it is simply a linear combination of the standard polynomial smoothness indicators $I[P_k]$, and it enters in the definition of the $\alpha_k$'s in order to drive the nonlinear weights $\omega_k$ closer to their optimal value $d_k$ in case of smoothness. In order to achieve its goals, $\tau$ should be at most of size $\mathcal{O}(h^{r+2})$ and explicit formulas for it are given in [CCD11]. In [DB13] it was shown that WENOZ can employ much smaller values for $\varepsilon$ and numerical tests show that this is beneficial for the error.

*2.1. CWENOZ reconstruction procedure*

Here we propose the introduction of the WENOZ weights in the CWENO reconstruction procedure.

**Definition 2** (CWENOZ). *For $r > 1$, consider the stencils $\mathbb{S}_0$ and the substencils $\mathbb{S}_1, \cdots, \mathbb{S}_r$ defined above, the polynomial $P_{\text{opt}}$ interpolating all the data in $\mathbb{S}_0$ and the polynomials $P_k$ interpolating the data in $\mathbb{S}_k$. Let also $d_0, d_1, \cdots, d_r$ be a set of positive, real coefficients $d_k \in (0,1)$, s.t. $\sum d_k = 1$. The CWENO operator computes a reconstruction polynomial*

$$P_{\text{rec}} = \mathsf{CWENOZ}(P_{\text{opt}}, P_1, \ldots, P_r) \in \mathbb{P}^{2r-2}$$

*as follows:*

1. *first, introduce the polynomial $P_0$ defined as*

$$P_0(x) = \frac{1}{d_0}\left(P_{\text{opt}}(x) - \sum_{k=1}^{r} d_k P_k(x)\right) \in \mathbb{P}^{2r-2}; \quad (6)$$

2. *then compute the nonlinear coefficients $\omega_k$ from the linear ones $d_k$ as*

$$\alpha_k^Z = d_k\left(1 + \left(\frac{\tau}{(I[P_k] + \epsilon)}\right)^t\right), \qquad \omega_k^Z = \frac{\alpha_k^Z}{\sum_{i=0}^{r} \alpha_i^Z}, \quad (7)$$

*where $I[P_k]$ denotes a regularity indicator (e.g. the Jiang-Shu ones of eq. (9)) evaluated on the polynomial $P_k$, $\tau$ is a suitable linear combination of $I[P_1], \ldots, I[P_r]$, $\epsilon$ is a small positive quantity and $t \geq 1$;*

3. *and finally define*

$$P_{\text{rec}}(x) = \sum_{k=0}^{r} \omega_k^Z P_k(x) \in \mathbb{P}^{2r-2}. \quad (8)$$

**Remark 2.** This new reconstruction technique shares the same advantages of CWENO listed in Remark 1. Additionally, the new definition of $\alpha_k^Z$ replacing the $\alpha_k$ of CWENO ensures more favourable convergence properties, as we sall prove below.

In this paper we consider as regularity indicators the classical ones proposed by Jiang-Shu in [JS96]

$$I[P] = \sum_{l \geq 1} h^{2l-1} \int_{\Omega_0} \left(\frac{\mathrm{d}^l}{\mathrm{d}x^l} P(x)\right)^2 \mathrm{d}x. \quad (9)$$

In view of the next results, let us introduce the notation $\theta(g(h))$ for the power of $h$ in the leading term of the Taylor series expansion of $g(h)$.

[DB13] propose two different definitions for $\tau$, that we adopt also for CWENOZ. The first one depends only on the parity of $r$, while the second one is taylored

| $r$ | $\tau$ | $\theta(\tau)$ | $\tau^{\text{opt}}$ | $\theta(\tau^{\text{opt}})$ |
|---|---|---|---|---|
| 2 | | | $\|I_1 - I_2\|$ | 3 |
| 3 | $\|I_1 - I_3\|$ | 5 | $\|I_1 - I_3\|$ | 5 |
| 4 | $\|I_1 - I_2 - I_3 + I_4\|$ | 6 | $\|I_1 + 3I_2 - 3I_3 - I_4\|$ | 7 |
| 5 | $\|I_1 - I_5\|$ | 7 | $\|I_1 + 2I_2 - 6I_3 + 2I_4 + I_5\|$ | 8 |
| 6 | $\|I_1 - I_2 - I_5 + I_6\|$ | 8 | $\|I_1 + I_2 - 8I_3 + 8I_4 - I_5 - I_6\|$ | 9 |

Table 1: The global smoothness indicator $\tau$ and the global optimal order smoothness indicator $\tau^{\text{opt}}$ and their leading truncation order $\theta(\tau), \theta(\tau^{\text{opt}})$ of the $(2r-1)$ order CWENOZ scheme. In the table we have used the notation $I_k$ for $I[P_k]$.

to each specific value of $r$ and aims at obtaining a quantity $\tau^{\text{opt}}$ which is smaller than the standard definition and yields advantages in the choice for $\varepsilon$ and in the global errors. $\tau^{\text{opt}}$ also allows to definition of CWENOZ 3, which is impossible since no $\tau$ would be zero for $r = 2$, according to the standard definition. In Table 1 we summarize the formulas of the $\tau$ and $\tau^{\text{opt}}$ with the respective leading truncation order $\theta(\tau)$.

The proof of the convergence of CWENOZ relies on the following result, that is proved in [DB13] for the finite differences setting. Let us focus on a reference cell $j = 0$ and assume that its cell centre is at $x_0 = 0$.

**Proposition 1 (cfr. Theorem 1 of [DB13]).** *Let $P$ be a polynomial of degree $q$ that interpolates the cell averages $\overline{u}_i$ of a function $u(x)$ in a stencil of $q+1$ continguous cells that include $\Omega_0$. The Jiang-Shu smoothness indicator (9) of $P$, can be written in bilinear form as*

$$I[P] = \langle \vec{w}, C\vec{w} \rangle$$

*where $C$ is a $q \times q$ positive semidefinite symmetric matrix and $\vec{w}$ is a vector such that $w_i = u^{(i)}(0)h^i + \mathcal{O}(h^{q+1})$, $i = 1, 2, \cdots, q$.*

PROOF. The positive semi-definiteness is trivial from the definition of the Jiang-Shu indicators.

The Jang-Shu indicator (9) of a generic polynomial of degree $q$, centered in 0, can be written in a bilinear form as

$$I\left[\sum_{l=0}^{q} a_l x^l\right] = \langle \vec{w}, C\vec{w} \rangle = \sum_{i=1}^{q}\sum_{j=1}^{q} C_{i,j} w_i w_j$$

where $C$ is a $q \times q$ positive semidefinite symmetric matrix, $\vec{w} \in \mathbb{R}^q$ is a vector whose components are given by $w_i = i! a_i h^i$, $i = 1, \cdots, q$. The entries of the upper part of $C = (C)_{i,j}$, $i = 1, \cdots, q$, $j = i, \cdots, q$ take the following expression

$$C_{i,j} = \begin{cases} \sum_{m=1}^{i} \frac{2^{2m-i-j}}{(i-m)!(j-m)!(i+j-2m+1)} & \text{if } i+j \text{ even} \\ 0 & \text{if } i+j \text{ odd.} \end{cases}$$

If the polynomial $\sum_{i=0}^{q} a_i x^i$ is interpolating the cell averages of a smooth enough function $u(x)$, see [CPSV16], its coefficients satisfy

$$a_i = \frac{1}{i!} u^{(i)}(0) + \mathcal{O}(h^{q-i+1}) \quad i = 0, 1, \cdots, q.$$

Then, in particular, we have the following expression for the components of the vector $\vec{w}$

$$w_i = u^{(i)}(0)\, h^i + \mathcal{O}(h^{q+1}) \quad i = 1, 2, \cdots, q.$$

**Proposition 2 (cfr. Theorem 7 of [DB13]).** *If $\theta(\epsilon) \leq \theta(\tau) - (r-1)/t$, then the CWENOZ scheme achieves the optimal order $2r-1$, as $h \to 0$, regardless of the presence of critical points.*

PROOF. The proof of Theorem 7 of [DB13], relies only on Theorem 1 of [DB13]. Since Proposition 1 ensures that the conclusions of Theorem 1 of [DB13] holds true also in the present finite volume setting, the proof of Proposition 2 can be easily obtained along the same lines as the proof of Theorem 7 of [DB13].

## 3. Spectral properties

The study of the performance of a numerical method starts from the behavior of the scheme on simple equations for which the exact solution is known. Since we are concentrating on the behavior of reconstruction algorithms, we consider only the effects of the reconstruction on the space approximation of differential operators.

As usual, we start from the simple linear advection equation $u_t + au_x = 0$, with periodic boundary conditions on $[0, 2\pi]$. The solution is found through Fourier series, as

$$u(x,t) = \sum_{-\infty}^{\infty} \hat{u}_k(t) e^{ikx} = \sum_{-\infty}^{\infty} \hat{u}_k(0) e^{ik(x-at)}.$$

The solution is obtained considering separately each Fourier mode, say $v_k(x,t) = \hat{u}_k(t) e^{ikx}$, substituting in the linear advection equation,

$$\frac{\mathrm{d}\hat{u}_k}{\mathrm{d}t} e^{ikx} + a\, ik\, \hat{u}_k e^{ikx} = 0,$$

and computing the evolution of the amplitude,

$$\hat{u}_k(t) = \hat{u}_k(0) e^{-ikat}.$$

If we consider instead a numerical scheme, based on a reconstruction in space, and an approximate derivative, the evolution of the amplitude will be given by

$$\frac{\mathrm{d}\hat{u}_k}{\mathrm{d}t} e^{ikx} + a\, \hat{u}_k D_x\left(e^{ikx}\right) = 0,$$

where $D_x$ is the discrete space derivative. If the scheme is based on a formula with constant coefficients, on the stencil $[x - mh, x - (m-1)h, \ldots, x + mh]$, we obtain

$$\frac{\mathrm{d}\hat{u}_k}{\mathrm{d}t} e^{ikx} + a\, \hat{u}_k \left( \sum_{\ell=-m}^{m} c_\ell e^{ikh\ell} \right) e^{ikx} = 0.$$

In other words, the function $e^{ikx}$ is an eigenfunction also for the discrete operator. Let $\tilde{\omega}_k$ be defined by the equation

$$\tilde{\omega}_k = \left( \sum_{\ell=-m}^{m} c_\ell e^{ikh\ell} \right) - ik. \tag{10}$$

Thus the amplitude of the Fourier mode, processed by the linear advection equation with the discrete derivative, is

$$\hat{u}_k(t) = \hat{u}_k(0) e^{-ikat} e^{-a\tilde{\omega}_k t}. \tag{11}$$

In other words, the effect of the numerical scheme is to modify the exact propagation of the Fourier mode, with the factor $e^{-a\tilde{\omega}_k t}$. More in detail, the exact propagation of the modified mode will be

$$v_k(x,t) = \hat{u}_k(0) e^{ik(x-\tilde{a}t)} e^{-a\mathrm{Re}(\tilde{\omega}_k)t}, \qquad \tilde{a} = a + \frac{a}{k}\mathrm{Im}(\tilde{\omega}_k). \tag{12}$$

Thus, the real part of $\tilde{\omega}_k$ controls the damping of the amplitude of the Fourier mode (provided that the scheme is stable). This is connected with the artificial diffusion of the scheme. The imaginary part of $\tilde{\omega}_k$ instead affects the propagation speed, as specified in (12), and this spurious behaviour determines the dispersive effects at the basis of the growth of spurious oscillations in non smooth solutions.

For example, for the first order upwind method, in the case of positive $a$, one gets

$$\begin{aligned}\tilde{\omega}_k &= \tfrac{1}{h}\left(1 - e^{-ikh}\right) - ik \\ &= -\tfrac{1}{2}k^2 h - i\tfrac{1}{6}k^3 h^2 + O(h^3),\end{aligned}$$

which shows that, as is well known, the main effect of this scheme is damping of high frequency modes. On the other hand, for the second order central derivative, one finds

$$\begin{aligned}\tilde{\omega}_k &= \tfrac{1}{h}\left(\frac{e^{ikh} - e^{-ikh}}{2}\right) - ik \\ &= -i\tfrac{1}{6}k^3 h^2 + O(h^4).\end{aligned}$$

Here, the spurious effect is an error in the propagation speed, given in (12). This of course is also well known, and can be found in several textbooks, such as [LV04] or [Tor09].

3.1. Diffusion and dispersion

Non linear schemes are more difficult to study, because naturally one cannot use the superposition principle. However, numerical diffusion and dispersion can be studied *a posteriori*, looking at how a given numerical scheme modifies single Fourier modes. This approach has been extensively considered in [Pir06], with the notion of ADR (Approximate Dispersion Relation). Here, we extend this study to include several variants of WENO schemes. Moreover, we develop a notion of *distortion*, to measure the non linear distorsive effects produced by high order non linear schemes.

Suppose we divide the interval $[-1, 1]$ in $2N + 1$ equal cells. We consider a semidiscrete scheme, in which the approximate derivative $D_x$ is computed with a non linear scheme. The single Fourier mode $\hat{u}_k(t)e^{ikx}$, $k = -N, \cdots, N$, is no longer an eigenfunction of $D_x$. We can however compute the Discrete Fourier Transform (DFT) of the output of the scheme, to see which Fourier modes have developed from the application of the non linear differentiation. Since the scheme works on real functions, compute, for $k = 1, \ldots, N$,

$$D_x \left[ \begin{array}{c} \sin(2\pi k x) \\ \cos(2\pi k x) \end{array} \right] = \sum_{\ell=1}^{N} \left[ \begin{array}{cc} \omega_{2\ell,2k} & \omega_{2\ell,2k+1} \\ \omega_{2\ell+1,2k} & \omega_{2\ell+1,2k+1} \end{array} \right] \left[ \begin{array}{c} \sin(2\pi l x) \\ \cos(2\pi l x) \end{array} \right]. \qquad (13)$$

Let $\Omega$ be the matrix with elements $\omega_{ij}$. The first column and the first row of $\Omega$ correspond to the constant mode, and thus they will be ignored in what follows. For a linear scheme only the 2 by 2 blocks along the diagonal of $\Omega$ would be non zero. Moreover, the exact derivative is the $2N$ by $2N$ block diagonal matrix

$$\mathbb{D} = \mathrm{diag}\,(E_k) = \mathrm{diag}\left( 2\pi k \left[ \begin{array}{cc} 0 & 1 \\ -1 & 0 \end{array} \right] \right), \qquad k = 1, \ldots, N.$$

It follows that all terms in $\Omega$ off the two main diagonals are spurious distortive effects. The error is given by the matrix $\Omega - \mathbb{D}$. We define the *relative* error due the non linear derivative as

$$E = |\Omega - \mathbb{D}|\,\mathrm{diag}\left( \tfrac{1}{2\pi k} \left[ \begin{array}{cc} 1 & 0 \\ 0 & 1 \end{array} \right] \right). \qquad (14)$$

Note that in the error matrix $E$, the error on each mode is normalized with its frequency, so that the elements of $E$ represent the relative errors on each mode.

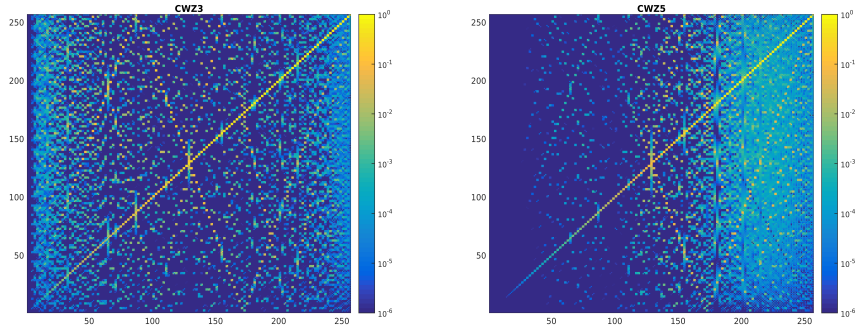

Figure 1: Error matrices $E$ for CWZ3 (left) and CWZ5 (right).

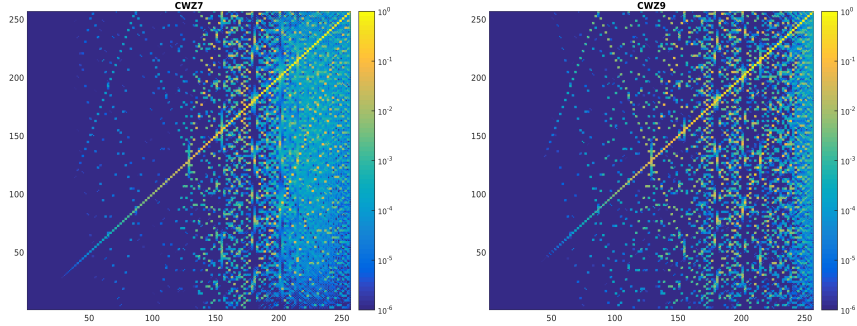

Figure 2: Error matrices $E$ for CWZ7 (left) and CWZ9 (right).

The two figures Fig. 1 and 2 show the size of the entries of the error matrices $E$ for low order and high order schemes, respectively, for $N = 128$, so that the total number of modes is 257. In particular, we exhibit the results obtained with the CWENOZ schemes of order 3 and 5 in Fig. 1, and the results yielded by CWENOZ 7 and 9 in Fig. 2. The errors on the low frequency modes appear in the left columns: as the order increases, more and more low frequency modes are well resolved. Note that in all cases the very high frequency modes are lost, and in fact the columns to the right of the matrix are filled with numerical artifacts. As a rule of thumb, you need at least 5 points per mode to resolve each wave, but the figures show that this rule can be relaxed for very high order schemes.

Along the two by two diagonal blocks, we find the errors on the $k$-th mode, which are the ones used to evaluate diffusion and dispersion errors, see for instance [Pir06]. However, the plots show clearly that we have not only diffusion and dispersion, but also the growth of spurious modes, which are a purely non linear effect. For low order modes, the error matrices $E$ are almost diagonal: this means that for low frequency modes distortion effects are minimal. As the order is increased, the number of modes which are distortion free increases. This is a further beneficial effect of high order schemes.

To recover diffusion and dispersion, we need to compute the derivative in complex form. Let $T$ be the matrix mapping the Fourier basis in complex form to the Fourier basis in real form. Then the matrix containing the coefficients of the approximate derivative $D_x$ on the basis $e^{ikx}$ can be written as $\Omega_\mathbb{C} = T^H \, \Omega \, T$. Let $\tilde{\omega}_{kj} = (\Omega_\mathbb{C})_{kj}$. Then,

$$\begin{aligned} \mathrm{Re}(\tilde{\omega}_{kk}) &\implies \text{diffusion error} \\ \mathrm{Im}(\tilde{\omega}_{kk}) - k &\implies \text{dispersion error}. \end{aligned}$$

We compare diffusion and dispersion for the CWENOZ schemes with standard WENO schemes, and with the compact WENO scheme (CWENO) of [CPSV16] (see the previous section). To minimize the number of possible de-grees of freedom, given by the possibility of choosing different parameters for different schemes, we always choose $\varepsilon = h^2$. In all cases, the abscissa carries the normalized Fourier mode, i.e. $\pi k/N$, where $N$ is the number of Fourier modes considered. The data on diffusion appear in Fig. 3 for third and fifth order schemes and in Fig. 4, for the schemes of order seven and nine. Diffusion is a purely numerical artifact. In all cases, the CWENOZ schemes have a slightly smaller amount of artificial diffusion than the other schemes, except on the very high frequency modes which in any case cannot be resolved on the chosen grid.

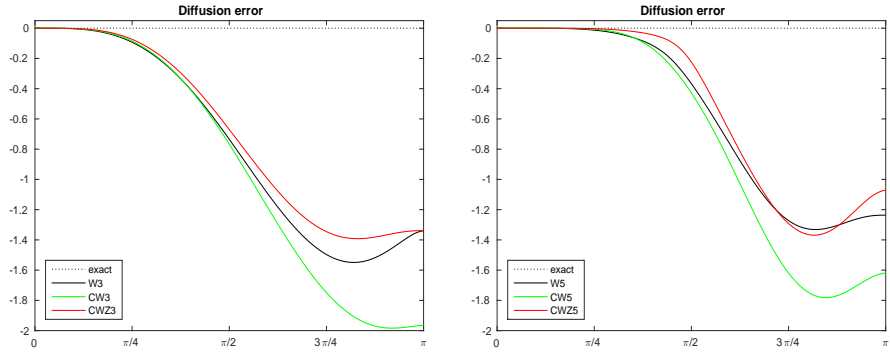

Figure 3: Diffusion error for schemes of order 3 (left) and 5 (right).

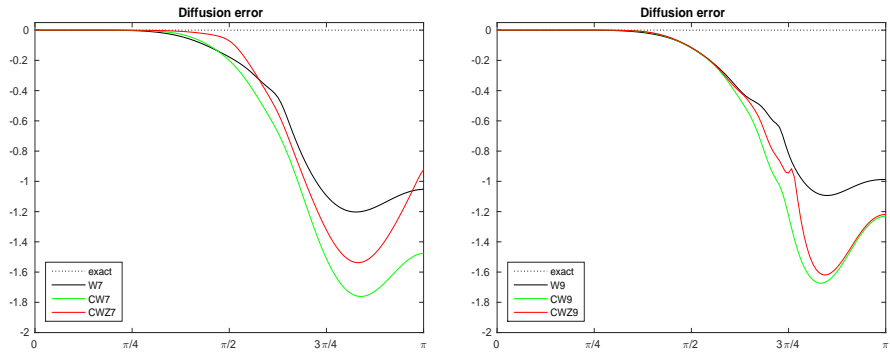

Figure 4: Diffusion error for schemes of order 7 (left) and 9 (right).

Fig. 5 shows the dispersive effects for the 3rd and 5th order schemes, while Fig. 6 refers to the 7th and 9th methods. In this case, the curves should be as close as possible to the dashed line, which contains the correct speed. It is clear that as the order increases, the numerical curves detach from the exact line at higher values of $k$. Here, both CWENO and CWENOZ seem to have an edge over standard WENO, especially in the high order case.

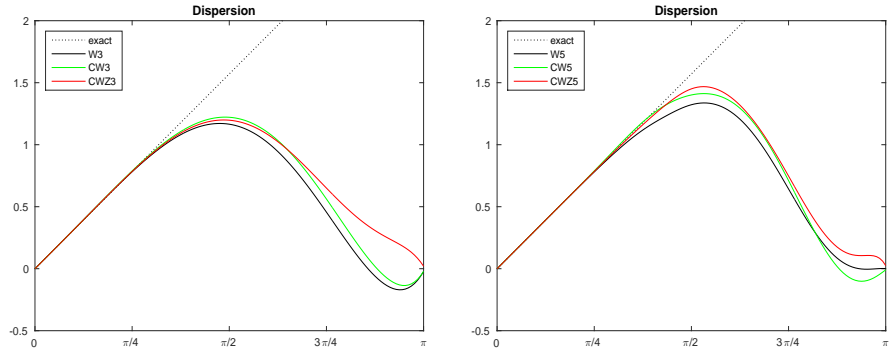

Figure 5: Dispersion for schemes of order 3 (left) and 5 (right).

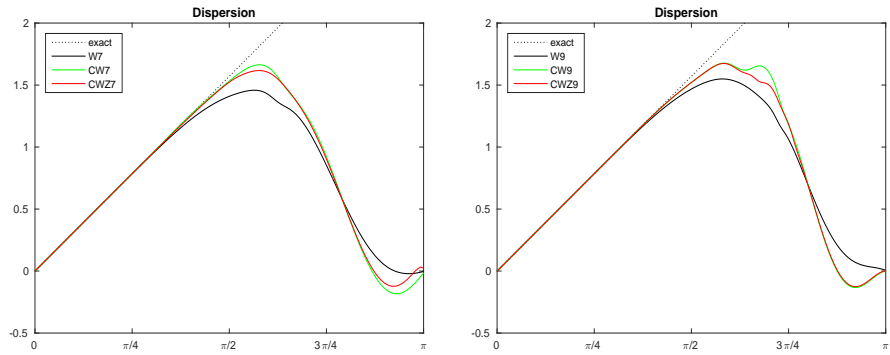

Figure 6: Dispersion for schemes of order 7 (left) and 9 (right).

### 3.2. Distortion and temperature

Diffusion and dispersion of numerical schemes have been studied by several authors, see for instance [Tre82], and are by now classical concepts, see [LV04] and references therein. They were derived for linear schemes, with the aid of the modified equation, and only later they were extended to non linear schemes, see [Pir06], or more recently [JGD15] or [BGS16]. To the best of our knowledge however, the non linear effects peculiar to non linear schemes such as WENO have not yet received the same consideration. However the error matrices $E$ of (14) and their plots in Figg. 1 and 2 clearly show that a large component of the error introduced by a scheme on a given mode consists in the development of spurious modes, which are not accounted for by the classical theory of diffusion and dispersion.

To account for this fact, we introduce two concepts which measure the distortion error on each mode and a global measure of distortion, which defines a parameter characterizing each scheme.

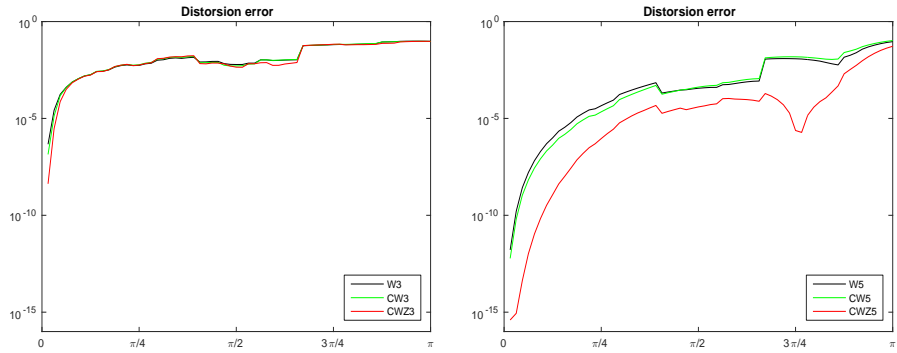

Figure 7: Distorsion for schemes of order 3 (left) and 5 (right).

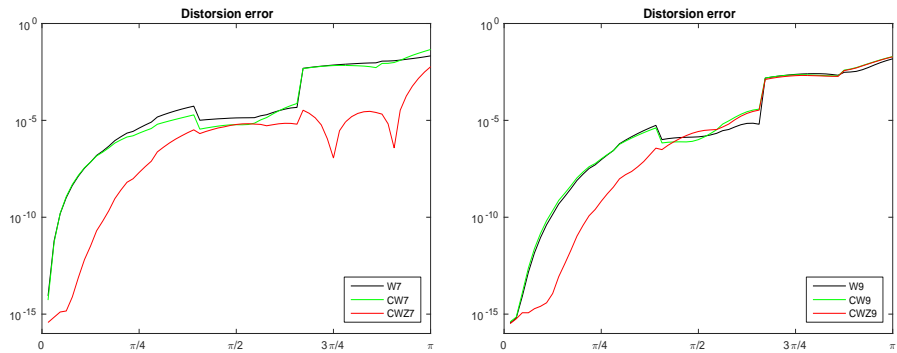

Figure 8: Distorsion for schemes of order 7 (left) and 9 (right).

**Definition 3 (Distortion error).** *The distortion error of a numerical scheme on the k-th mode of the Fourier basis is defined as*

$$\delta_k = \frac{1}{N} \sum_{\ell \neq k} |(\Omega_\mathbb{C})_{\ell k}|, \tag{15}$$

*where $\Omega_\mathbb{C}$ is the matrix containing the numerical derivative in the Fourier basis, defined in the previous subsection.*

The scaling introduced in the definition ensures that the plots of the distortion error do not depend on the number of nodes of the grid. Thus these figures are characteristic of each scheme: they do not depend on the grid spacing, but only on the degree of the interpolating polynomial, or the choice of parameters in the non linear weights. They are a sort of signature of the scheme. The distortion errors of WENO type schemes of a given order are compared in Figg. 7 for low order and 8 for high order. We see that we start with very small distortion errors on low frequency modes, which increase with a power law as the wave number is increased. Note that the distortion error again is smaller for high order schemes, on low frequency modes. More in detail: if we fix a certain threshold, say for instance Tol, the frequency $K$ for which the distortion error is $\delta_k <$ Tol is an increasing function of the accuracy $p$ of the scheme.

We also stress that the CWENOZ schemes are by far less distorsive than CWENO and WENO schemes, except at third order, where they almost coincide.

The distorsion error suggests the introduction of a single parameter, measuring the non linear distorsive effects of a numerical scheme.

**Definition 4 (Temperature).** *The temperature of a numerical scheme on the $k-$th mode is*

$$T_k = \frac{1}{N^3} \sum_{\ell=1}^{N} \tilde{\omega}_{k\ell}((k-\ell)/\pi)^2.$$

*The j-th temperature of the scheme is*

$$\mathcal{T}_j = \frac{1}{j} \sum_{\ell=1}^{j} T_\ell,$$

*i.e. it is the average temperature of the first j modes. Finally the temperature $\mathcal{T}$ of a scheme is defined as*

$$\mathcal{T} = \mathcal{T}_{N/2} \tag{16}$$

Note that the temperature $\mathcal{T}$ is defined using only the modes that can be resolved on a given grid. Moreover, with these scalings, the temperatures $T_k$ and $\mathcal{T}$ depend only on the particular scheme chosen, and not on the number of grid points. The definition of $\mathcal{T}$ yields a single parameter which characterizes the distortion errors of a scheme. For a linear scheme or for the exact derivative, $\mathcal{T}_j = 0, \forall j$ and in particular $\mathcal{T} = 0$.

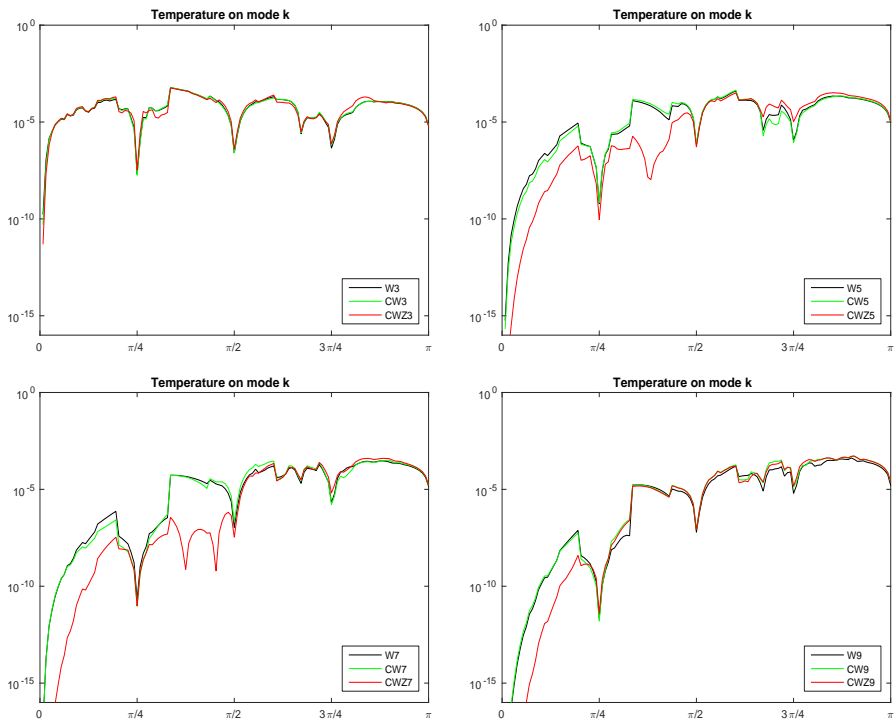

Figure 9: Temperature $T_k$ for several schemes: order 3 (left) and 5 (right) in the top row, order 7 (left) and 9 (right) in the bottom row.

|        | 3       | 5       | 7       | 9        |
|--------|---------|---------|---------|----------|
| WENO   | 4.60e-5 | 1.14e-6 | 9.31e-8 | 7.17e-9  |
| CWENO  | 5.06e-5 | 7.56e-7 | 3.52e-8 | 5.61e-9  |
| CWENOZ | 5.37e-5 | 7.30e-8 | 4.13e-9 | 4.11e-10 |

Table 2: Temperature $\mathcal{T}$ of WENO schemes, for orders 3 to 9.

Fig. 9 contains the temperatures $T_k$ for each mode of all schemes studied here. The notion of temperature weighs more the spurious modes which are far away in frequency space from the exact mode. More precisely, it measures the variance of the numerical derivative around the exact mode, and in this sense it reminds of the classical notion of temperature.

Linear schemes, which present no distortion error, have zero temperature, but are oscillatory. Thus a good scheme is "cool", in the sense that its temperature must be as low as possible, without developing spurious oscillations. From the plots, we see that CWENOZ has the smallest temperature among the essentially non-oscillatory schemes considered in this work. This is due to the fact that the weights in CWENOZ tend to privilege the choice of the high order polynomial entering in the reconstruction, thus approaching a linear scheme. Note that the temperature decreases, as the order is increased.

Table 3.2 contains the temperature $\mathcal{T}$ of (16). The normalization appearing in the definition of temperature ensures that these numbers characterize a numerical scheme, because $\mathcal{T}$ does not depend on the number of grid nodes on which the method is tested. The temperature of the various schemes tested decrease with the order of the scheme. In particular, CWENOZ is the coolest among all schemes tested, at each order.

All plots shown are obtained with $N = 128$ real modes (plus the zero mode), that is the number of grid nodes is $2N + 1 = 257$, but the plots obtained with different values of $N$ can be superposed almost exactly, thanks to the normalization chosen.

In all tests, we used the same choice of the parameter $\varepsilon$ in all schemes. Modifying the choice of $\varepsilon$ does not produce significant differences.

### 4. Numerical tests

In this section we perform a number of tests on scalar and systems of conservation and balance laws, in order to compare schemes based on the reconstructions studied in this paper, corroborating the findings of the previous sections.

Computational grids are set up subdividing the domain in $M$ uniform cells and considering, when appropriate, ghost cells outside the physical domain, whose cell averages are extrapolated from the cell averages of the inner cells, taking into account the boundary conditions. Since the purpose of the tests is to compare the reconstructions, all other ingredients of the numerical scheme were chosen as simple as possible and the same for all reconstructions. In particular the numerical flux is the Local Lax Friedrichs, the CFL number was always set to 0.45 and the ODE system obtained by semidiscretization was solved numerically with a Runge-Kutta scheme of order matching the order of accuracy of the spatial reconstruction. More precisely, the third order scheme employs the classical third order (strong stability preserving) SSP Runge-Kutta with three stages [JS96], the fifth order scheme the fifth order Runge-Kutta with six stages of [But08, §3.2.5], the scheme of order seven relies on the nine-stages Runge-Kutta of [But08, pag 196] and the scheme of order nine employs the Runge-Kutta with eighteen stages of order ten of [Cur75]. Clearly, other Runge-Kutta or multistep schemes and different Riemann solvers could be used instead.

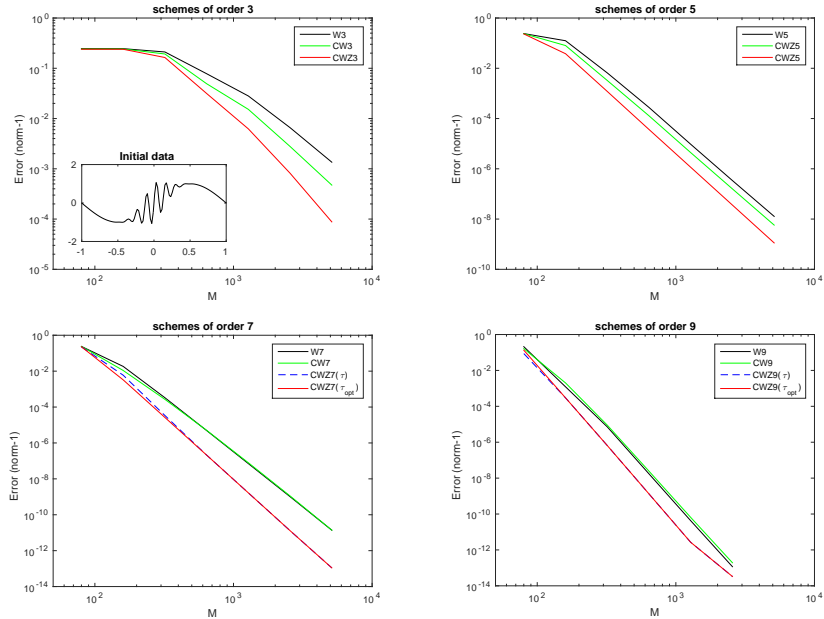

Figure 10: Convergence test for the linear transport of the smooth data (17) for schemes of order 3, 5, 7 and 9. The initial data is shown in the inset of the top-left panel.

### 4.1. Linear transport

We start with a set of tests on the linear transport equation $u_t + u_x = 0$ on the domain $[-1, 1]$ with periodic boundary conditions.

We first consider the initial data

$$u(x) = \sin(\pi x) - \sin(15\pi x)e^{-20x^2} \qquad (17)$$

already introduced in [SCR16], that is a signal that mixes low and high frequency components and contains all Fourier modes. Fig. 10 shows the convergence of the errors for the schemes of order 3 to 9. In all tests, one obviously observes that the error does not decrease until the grid is fine enough to resolve details in the middle of the initial data. Thus the convergence with the expected rate starts approximately for $M > 320$ for the third order schemes, for $M > 160$ for fifth order ones and for $M > 80$ at seventh and ninth order. The graphs clearly show that the CWENOZ schemes outperform the other reconstructions because they provide the smallest errors and the best convergence rates, since the optimal rates appear already for relatively small $M$. CWENO 3 and CWENO 5 perform better that the WENO schemes of corresponding order, while CWENO 7 and CWENO 9 are on par with the corresponding WENO schemes. Note finally that the choice of $\tau^{\text{opt}}$ in CWENOZ 7 and CWENOZ 9 gives an advantage for small $M$ and this is important since in practice one would like to employ high order schemes especially on coarse grids. For CWENOZ 9 the test stops at $M = 2560$ instead of 5120 because at this resolution the error is already close to the machine precision.

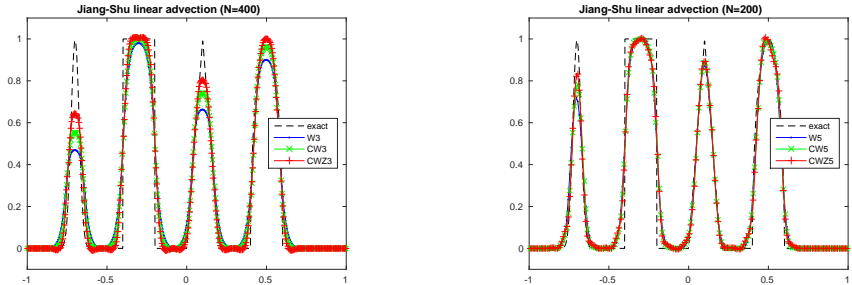

Figure 11: Linear advection of a smooth bump, square wave, triangular wave and a semi-ellipse. Left: schemes of order 3, with 400 points. Right: schemes of order 5, with 200 points.

*Non-smooth initial data.* Next, we consider the initial data proposed in [JS96], which contain a smooth part, a square wave, a triangular wave and a semi-ellipse in the periodic domain $[-1, 1]$. The solutions obtained at time $t = 8$ with the order 3 and order 5 schemes are plotted in Fig. 11. 400 points have been used in the first case and 200 in the second one, in order to be able to distinguish the curves. One can see that the CWENOZ schemes resolve much better the smooth part of the data (Gaussian bump on the left) and the amplitude of the triangular wave (third bump), without oscillating on the jumps of the square wave and close to the steep gradient of the semi-ellipse.

In Fig. 12 we plot the nonlinear weights employed by the reconstructions on the initial data. In particular we consider the relative error $(\omega_* - d_*)/d_*$ for the central $P_0$ for CWENO and CWENOZ schemes and for the left-most polynomial for WENO schemes, which do not have necessarily a central polynomial. The vertical scale is logarithmic and its minimum is set to the machine precision; missing symbols are exact zero values. One can note that, especially in the middle of the Gaussian bump and of the semi-ellipse, CWENOZ is much more effective than CWENO in using the optimal weights for the reconstruction.

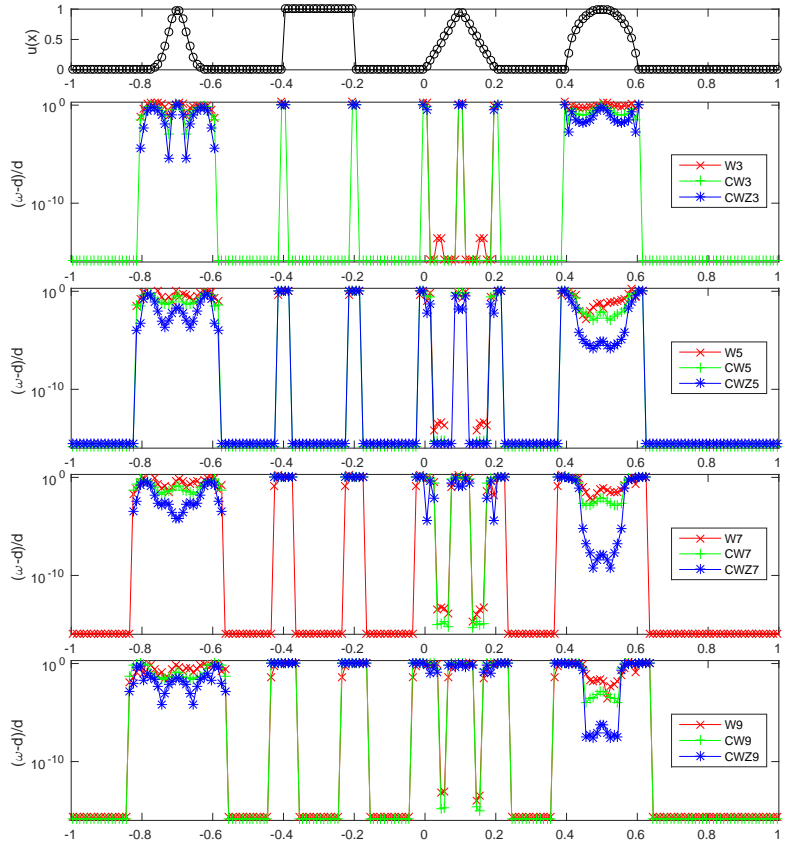

Figure 12: Top panel: initial data. Lower panels: relative errors on the nonlinear weights for schemes of order 3 to 9, from top to bottom. We plot $(\omega_1 - d_1)/d_1$ for WENO schemes and $(\omega_0 - d_0)/d_0$ for CWENO and CWENOZ schemes. (Missing symbols correspond to exact zero values.)

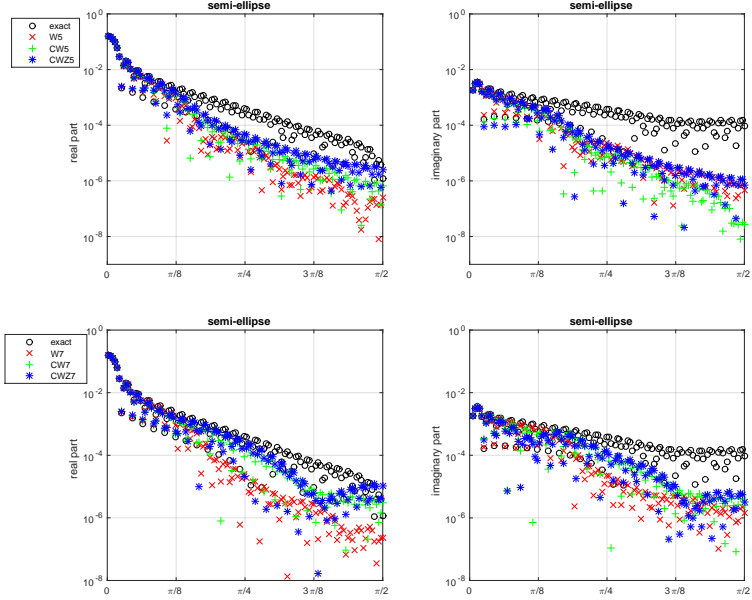

Figure 13: Real (left) and imaginary (right) part of the DFT of the exact and of the numerical solution of linear advection with data (18).

These results corroborate the findings of the study of the distorsion (Figg. 7, 8) and of the temperature (Figg. 9) of the reconstruction procedures. The lower errors of the weights employed by the CWENOZ scheme clearly correspond to their lower distorsion and smaller temperature.

The effects of the discrepancies between the nonlinear weights and the optimal ones is further studied in Fig. 13, which shows the DFT of the exact and of the computed solutions on the semi-ellipse Jiang-Shu data of Fig. 12. Namely, we run the code with initial data

$$u_0(x) = 1/6(F(x,10,0.5-\delta) + F(x,10,0.5+\delta) + 4F(x,10,0.5)) \qquad (18)$$

with $F(x,\alpha,a) = \sqrt{\max(1-\alpha^2(x-a)^2,0)}$ and $\delta = 0.005$. For symmetry we show only half of the spectrum in Fig. 13 for schemes of order 5 and 7. Note that the data from CWENOZ schemes are the closest to the exact DFT even for high frequencies.

Finally, we study the distorsion of linear advection with initial data (17) already used in the convergence plots. Fig. 14 contains the DFT of the exact and of the numerical solutions of schemes of order 7 and 9, as well as the relative error between the nonlinear weights and the optimal weights. Again, we see that a bias towards the optimal weights obtained with CWENOZ results in a lower distorsion of the signal.

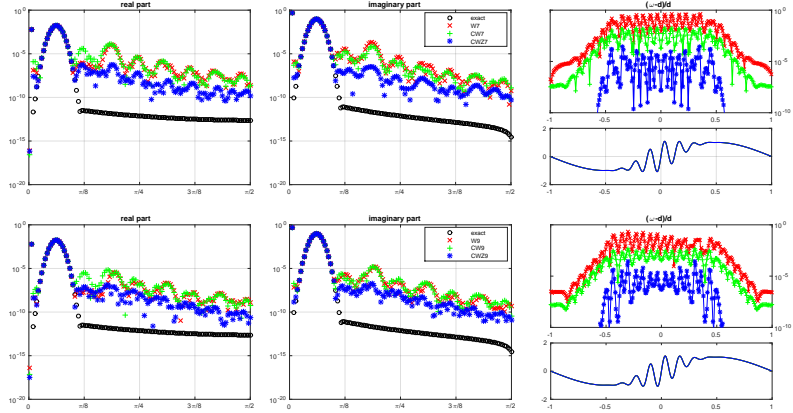

Figure 14: Comparison of schemes of order 7 and 9 on the solution of linear advection with initial data (17). Left and middle: real and imaginary part of Fourier transform of final data. Right: profile (bottom) and relative error on nonlinear coefficients (top).

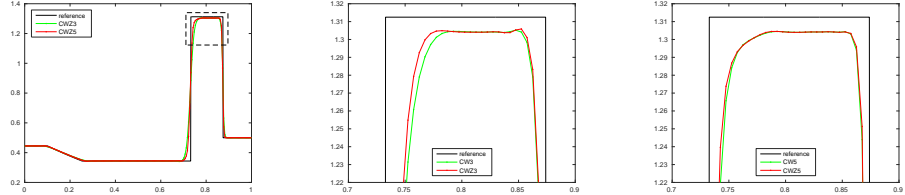

Figure 15: Lax shock tube with 200 grid points. Left: comparison of third and fifth order CWENOZ schemes. Middle and right: comparison of CWENO and CWENOZ schemes of order 3 and 5 on the detail highlighted by the dashed rectangle in the left panel.

4.2. Euler gas dynamics

Next we consider the Euler system of one-dimensional gas-dynamics, in order to study the non-oscillatory properties of the schemes on nonlinear problems. The systems of conservation laws is

$$\partial_t \begin{pmatrix} \rho \\ \rho u \\ E \end{pmatrix} + \partial_x \begin{pmatrix} \rho u \\ \rho u^2 + p \\ u(E + p) \end{pmatrix} = 0,$$

where $\rho$ is the gas density, $u$ the velocity, $p$ the pressure, and $E$ the energy per unit volume. The pressure is linked to the other variables through the equation of state of an ideal gas, namely $p = (E - \frac{1}{2}\rho u^2)(\gamma - 1)$, and we take $\gamma = 1.4$. As usual with very high order schemes we apply the reconstructions to the local characteristic variables.

*Lax's shock tube.* In order to demonstrate the essentially non-oscillatory properties of the proposed schemes, we consider the Riemann problem by Lax, which has the following left and right states: $\rho_L = 0.445, u_L = 0.6989, p_L = 3.5277$ and $\rho_R = 0.5, u_R = 0, p_R = 0.571$. The solution develops a left-moving rarefaction and two right-moving waves: a contact discontinuity and a shock, separated by a constant state. Capturing correctly and without spurious oscillations the plateaux between the contact and the shock wave is a hard test for the essentially non-oscillatory properties of a scheme. In Fig. 15 we compare the CWENOZ 3 and the CWENOZ 5 scheme, showing the increased resolution of the higher order scheme and the absence of spurious oscillations. The middle and right panels present details of the area in the dashed rectangle of the whole solution shown in the left panel.

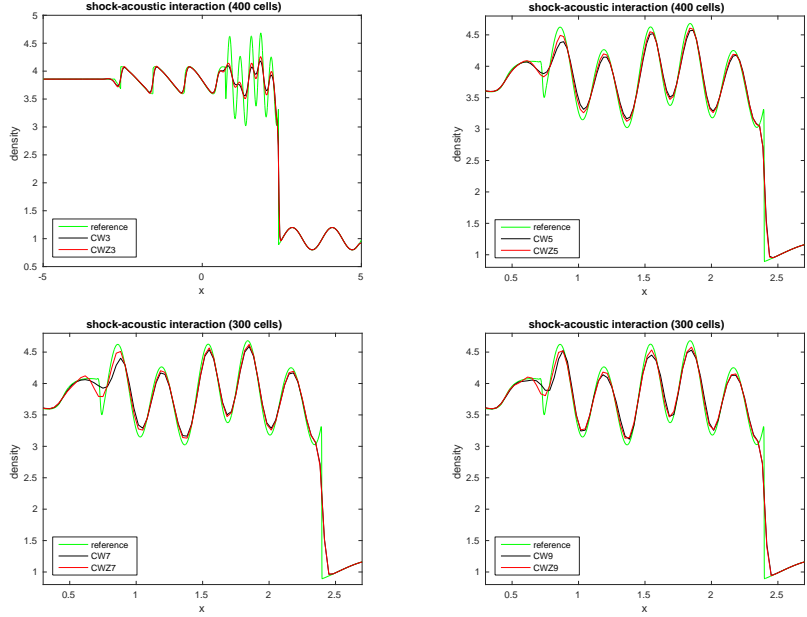

Figure 16: Density for the shock-acoustic interaction problem by CWENO and CWENOZ schemes of order 3 (top-left), 5 (top-right), 7 (bottom-left) and 9 (bottom-right). The top left panel shows the solution in the whole domain for the third order schemes. The other panels show only an enlargement of the turbulence area. The solution with the schemes of order 3 and 5 employed 400 cells, the solutions with the schemes of order 7 and 9 employed 300 cells.

*Interaction of a shock with a standing acoustic wave.* This test was proposed in [SO88] and consists in the following initial data

$$(\rho, v, p) = \begin{cases} (3.857143, 2.629369, 10.333333), & x \in [-5, -4) \\ (1.0 + 0.2\sin(5x), 0.0, 1.0), & x \in [-4, 5] \end{cases}$$

in the domain $[-5, 5]$, computing the evolution until $t = 1.8$. The solution is plotted in the top-left panel of Fig. 16, with a reference solution in light green and the solutions obtained with the CWENO 3 and CWENOZ 3 schemes in black and red respectively. The difficulty in this test is to avoid both spurious oscillations and excessive numerical diffusion at the shock, which allows to capture correctly the high frequency waves that form behind the strong shock, as a result of its interaction with the acoustic wave. In the remaining panels of Fig. 16 we thus show only a zoom of the solution in this problematic area. In order to highlight the differences among the schemes, the number of cells is considerably lower than the typical values used in this test. In particular, note the very tiny peak immediately behind the shock cannot be captured by any of the schemes due to the lack of grid resolution.

The top-left panel of Fig. 16 shows the solution obtained with the third order scheme. Here we see that CWENOZ 3 approximates slightly better the amplitude of the small waves in the region with $x \in [0.7, 2.3]$. Of course using more than 400 cells quickly allows the scheme to capture these waves correctly. The top-right panel shows the results for the order 5 schemes, again with 400 cells. It is apparent that the fifth order schemes can capture much more accurately the amplitude of the high frequency waves. Here again, CWENOZ 5 is more accurate than CWENO 5. The bottom panels are about the computations with the order 7 and 9 schemes. In order to distinguish the computational results from the reference solution, the grid resolution was lowered to a mere 300 cells. Despite this, the little waves are captured better than by the fifth order schemes and here again CWENOZ is slightly better than CWENO (note in particular the first wave on the left, at about $x = 0.75$).

### 4.3. Shallow water equations

Finally, since both CWENO and CWENOZ reconstructions are able to provide uniform accuracy also within the computational cells, we compare the reconstructions on a well-balanced scheme for the shallow water system, which makes use of quadrature points in the interior of the computational cells to integrate the source term. We thus consider the system

$$\partial_t \begin{pmatrix} h \\ q \end{pmatrix} + \partial_x \begin{pmatrix} q \\ q^2/h + \frac{1}{2}gh^2 \end{pmatrix} = \begin{pmatrix} 0 \\ -ghz_x \end{pmatrix},$$

where $h$ denotes the water height, $q$ is the discharge, $z(x)$ the known bottom topography and $g$ is the gravitational constant.

The well balanced quadrature is computed using Richardson's extrapolation, based on the trapezoidal rule [NPPN06]. This means that the source term average is computed using the two boundary value reconstructions and additionally 3, 7 and 15 internal nodes to achieve 5th, 7th and 9th order accuracy respectively. We emphasise that all these reconstructed data are computed from a single CWENO or CWENOZ reconstruction polynomial, using the same weights for all reconstruction points. On the other hand with WENO each node would require a separate reconstruction step.

In order to compare the accuracy of the schemes, we compute the flow proposed in [XS05], i.e. with initial data given by

$$z(x) = \sin^2(\pi x) \qquad h(0, x) = 5 + e^{\cos(2\pi x)} \qquad q(0, x) = \sin(\cos(2\pi x)), \qquad (19)$$

|      | CW3       |      | CWZ3      |      | CW5       |      | CWZ5      |      |
|------|-----------|------|-----------|------|-----------|------|-----------|------|
| N    | error     | rate | error     | rate | error     | rate | error     | rate |
| 32   | 9.37e-03  |      | 6.65e-03  |      | 4.01e-04  |      | 2.58e-04  |      |
| 64   | 1.44e-03  | 2.70 | 7.48e-04  | 3.15 | 1.73e-05  | 4.54 | 9.72e-06  | 4.73 |
| 128  | 1.56e-04  | 3.21 | 6.40e-05  | 3.55 | 5.74e-07  | 4.91 | 3.15e-07  | 4.95 |
| 256  | 1.57e-05  | 3.31 | 7.17e-06  | 3.16 | 1.81e-08  | 4.98 | 9.99e-09  | 4.98 |
| 512  | 1.83e-06  | 3.10 | 8.80e-07  | 3.03 | 5.70e-10  | 4.99 | 3.13e-10  | 5.00 |
| 1024 | 2.29e-07  | 3.00 | 1.10e-07  | 3.00 | 1.79e-11  | 5.00 | 9.80e-12  | 5.00 |
|      | CW7       |      | CWZ7      |      | CW9       |      | CWZ9      |      |
| N    | error     | rate | error     | rate | error     | rate | error     | rate |
| 16   | 1.30e-03  |      | 1.63e-03  |      | 7.02e-04  |      | 6.88e-04  |      |
| 32   | 7.25e-05  | 4.17 | 6.22e-05  | 4.71 | 2.82e-05  | 4.64 | 2.38e-05  | 4.85 |
| 64   | 6.70e-07  | 6.76 | 7.44e-07  | 6.39 | 1.22e-07  | 7.85 | 1.17e-07  | 7.67 |
| 128  | 5.02e-09  | 7.06 | 6.68e-09  | 6.80 | 3.44e-10  | 8.47 | 3.15e-10  | 8.54 |
| 256  | 3.91e-11  | 7.00 | 5.37e-11  | 6.96 | 7.43e-13  | 8.86 | 6.65e-13  | 8.89 |
| 512  | 3.07e-13  | 6.99 | 4.25e-13  | 6.98 | *7.11e-15* | *6.71* | *6.91e-15* | *6.59* |

Table 3: Errors in the SWE smooth test. Data in italic are most likely affected by machine precision.

with periodic boundary conditions on the domain $[0, 1]$. At time $t = 0.1$ the solution is still smooth and we compare the numerical results with a reference solution computed with a fourth order scheme on 32768 cells. The 1-norm of the errors appears in Table 3.

The upper part of Table 3 concerns the schemes of order 3 and 5. In both cases the CWENOZ errors are lower than the corresponding CWENO data. In the bottom part of the table, we see that the schemes of order 7 and 9 have the expected accuracy, but now the errors are very close. Clearly, as the error approaches machine precision, the convergence rate deteriorates. Note that for the high order schemes, the number of computational cells was lowered by a factor of 2, because the optimal order of accuracy is achieved faster.

## 5. Conclusions

In this paper we considered several high order essentially non-oscillatory reconstructions of WENO type for finite volume schemes applied to balance laws. In order to compare high order schemes, the concepts of numerical diffusion and numerical dispersion have been extended to the nonlinear case, to account for the high nonlinearity of such schemes. Unlike [Pir06], however, we compute the Approximate Dispersion Relations on the discrete derivative itself, to avoid polluting contributions from the time integrators.

The notions of numerical diffusion and dispersion generalize the idea of the von Neumann analysis. In the linear setting a Fourier mode is an eigenfunction both for the exact and the discrete derivative and only the eigenvalue is different. However, in the nonlinear case, there are distorsive effects: the dis-crete operator applied to a single Fourier mode creates new modes which are spurious effects. These obviously cannot be measured with tools derived from linear analysis. As far as we know, there is no study of such distorsive effects in the literature. For this reason, in this paper, we have introduced the notion of distorsion error, which is a measure of the amplitude of the spurious modes generated by the reconstruction, and a notion of temperature, which measures the distance in frequency space of the spurious modes from the exact signal. Both these quantities are scaled in such a way that they do not depend on the grid spacing and thus they provide measures of the distorsive effects which are characteristic of each scheme.

A linear scheme has zero temperature and zero distorsive errors, but it may be oscillatory. Nonlinear tools designed to prevent spurious oscillations increase the temperature of a scheme, which reads the distorsive effects. Thus, we seek for a scheme which is cold, but distorsive enough to prevent oscillations. For this reason we considered the WENOZ choice of weights, which is designed to reduce the nonlinearity of the original WENO schemes. We exported this idea to the new CWENO schemes introduced in [CPSV16], defining the CWENOZ reconstructions. In this fashion we enjoy the low temperature of the WENOZ weights and the possibility of creating uniformly accurate reconstructions of the CWENO approach.

Several tests illustrate the performance of the CWENOZ and the ability of the ideas of distorsion error and temperature to predict the behavior of the schemes.

*Acknowledgments.* This work was supported by "National Group for Scientific Computation (GNCS-INDAM)".